\documentclass[12pt]{article}
\usepackage{latexsym,cite,amssymb,amsmath,enumerate,geometry}
\newtheorem{theorem}{Theorem}

\newtheorem{lemma}[theorem]{Lemma}

\newtheorem{claim}{Claim}

\usepackage{tikz}

\begin{document}

\title{Largest Domination Number and Smallest Independence Number of Forests with given Degree Sequence}
\author{Michael Gentner$^1$, Michael A. Henning$^2$, and Dieter Rautenbach$^1$}
\date{}
\maketitle
\vspace{-10mm}
\begin{center}
{\small
$^1$ Institute of Optimization and Operations Research, Ulm University, Ulm, Germany\\
\texttt{michael.gentner@uni-ulm.de, dieter.rautenbach@uni-ulm.de}\\[3mm]
$^2$ Department of Mathematics, University of Johannesburg, Auckland Park, 2006, South Africa\\
\texttt{mahenning@uj.ac.za}}
\end{center}

\begin{abstract}
For a sequence $d$ of non-negative integers,
let ${\cal F}(d)$ be the set of all forests whose degree sequence is $d$.
We present closed formulas for 
$\gamma_{\max}^{\cal F}(d)=\max\{ \gamma(F):F\in {\cal F}(d)\}$ 
and
$\alpha_{\min}^{\cal F}(d)=\min\{ \alpha(F):F\in {\cal F}(d)\}$
where $\gamma(F)$ and $\alpha(F)$ are the domination number and the independence number of a forest $F$, respectively.
\end{abstract}

{\small 


}

\pagebreak

\section{Introduction}

We consider finite, simple, and undirected graphs, and use standard terminology.
For a sequence $d$ of non-negative integers, let ${\cal G}(d)$ be the set of all graphs with degree sequence $d$. Similarly, let ${\cal F}(d)$ be the set of all forests with degree sequence $d$.
For some graph parameter $\pi$ and an optimization goal ${\rm opt}\in \{\min,\max\}$, let
$$\pi_{\rm opt}(d)={\rm opt}\{\pi(G):G\in {\cal G}(d)\}\mbox{ and }
\pi_{\rm opt}^{\cal F}(d)={\rm opt}\{\pi(F):F\in {\cal F}(d)\}.$$
Note that for every graph $G$ with degree sequence $d$,
the values of $\pi_{\min}(d)$ and $\pi_{\max}(d)$
are the best possible lower and upper bounds on $\pi(G)$
that only depend on the degree sequence of $G$.

In the present paper we focus on two of the most prominent computationally hard graph parameters;
the domination number $\gamma(G)$ and the independence number $\alpha(G)$ of a graph $G$.
Many of the well known bounds \cite{m,c,fms,t,w,dhh,hhs,l,bhks,p,lp}
on these two parameters
depend only on the degree sequence,
or on derived quantities such as the order, the size, the minimum degree, and the maximum degree,
which motivates the study of $\pi_{\min}(d)$ and $\pi_{\max}(d)$.
Rao \cite{r} obtained the surprising result that $\alpha_{\max}(d)$ can be determined efficiently for every degree sequence $d$
(cf. also \cite{r2,kl,y}). 
In \cite{ghr} we showed that $\gamma_{\min}(d)$
can be determined efficiently for degree sequences with bounded entries,
and we gave closed formulas for $\gamma_{\min}^{\cal F}(d)$ as well as for  $\alpha_{\max}^{\cal F}(d)$.

Bauer et al. \cite{bhks} conjectured that $\alpha_{\min}(d)$ is computationally hard,
and we \cite{ghr} believe that the same is true for $\gamma_{\max}(d)$.
Therefore, for these last two parameters,
we focus on the more restricted case of forests.
Our main results are closed formulas for
$\gamma_{\max}^{\cal F}(d)$
and
$\alpha_{\min}^{\cal F}(d)$.
Note that for some degree sequences of forests,
there are exponentially many non-isomorphic realizations.
Therefore, the simple linear time algorithms
that determine the domination number and the independence number
of a given forest do not lead to an efficient algorithm that determines
$\gamma_{\max}^{\cal F}(d)$
and
$\alpha_{\min}^{\cal F}(d)$.

Let $d$ be a sequence $(d_1,\ldots,d_n)$ of $n$ non-negative integers.
The sequence $d$ is non-increasing if $d_1\geq d_2\geq \ldots \geq d_n$.
For a non-negative integer $i$, let $n_i(d)$ and $n_{\geq i}(d)$
be the numbers of entries of $d$ that are equal to $i$ and at least $i$, respectively.
It is well-known that $d$ is the degree sequence of some forest
if and only if $\sum_{i=1}^nd_i$ is an even number at most $2(n-n_0(d))-2$.
More specifically, if $\sum_{i=1}^nd_i=2(n-n_0(d))-2c$ for some positive integer $c$,
then every forest with degree sequence $d$ has
$n_0(d)$ isolated vertices and $c$ further non-trivial components.
In particular, if all entries of $d$ are positive, then $d$ is the degree sequence of a tree
if and only if $\sum_{i=1}^nd_i=2n-2$.

Let $G$ be a graph.
For a non-negative integer $i$,
let $V_i(G)$ and $V_{\geq i}(G)$ be the sets of vertices of $G$
of degree $i$ and at least $i$, respectively.
A vertex of degree at least $2$ with a neighbor of degree $1$ is a support vertex.
A dominating set of $G$ is a set $D$ of vertices of $G$
such that every vertex of $G$ that does not lie in $D$ has a neighbor in $D$,
and the domination number $\gamma(G)$ of $G$
is the minimum cardinality of a dominating set of $G$.
An independent set in $G$ is a set of pairwise non-adjacent vertices of $G$,
and the independence number $\alpha(G)$ of $G$
is the maximum cardinality of an independent set in $G$.

\section{Results}

We begin with two preparatory lemmas.

\begin{lemma}\label{lemma1}
If $T$ is a tree of order $n$,
then there is a set $D$ of at most $\left\lceil\frac{n-2}{3}\right\rceil$ vertices of $T$
such that every vertex $u$ of $T$
that has degree at least $2$ and does not belong to $D$
has a neighbor in $D$.
\end{lemma}
{\it Proof:}
We prove the statement by induction on the order $n$.
If $n\leq 2$,
then $T$ has no vertex of degree at least $2$,
and $D=\emptyset$ has the desired properties.
Now, let $n\geq 3$.
Let $u_0u_1\ldots u_{\ell}$ be a longest path in $T$.
If $\ell=2$, then $T$ is a star of order at least $3$
with a center vertex $u$, and $D=\{ u\}$ has the desired properties.
Hence, we may assume that $\ell\geq 3$.
Let $T'$ be the component of $T-u_2u_3$ that contains $u_3$.
Clearly, the order $n'$ of $T'$ satisfies $n'\leq n-3$.
By induction,
there is a set $D'$ of at most $\left\lceil\frac{n'-2}{3}\right\rceil$ vertices of $T'$
such that every vertex $u$ of $T'$
that has degree at least $2$ (in $T'$)
and does not belong to $D'$
has a neighbor in $D'$.
Now, the set $D=D'\cup \{ u_2\}$
contains at most
$\left\lceil\frac{n'-2}{3}\right\rceil+1\leq \left\lceil\frac{n-2}{3}\right\rceil$
vertices and has the desired properties.
Note that $u_3$ might have degree less than $2$ in $T'$
but is adjacent to $u_2\in D$ in $T$.
$\Box$

\begin{lemma}\label{lemma2}
If $d=(d_1,\ldots,d_n)$ is a non-increasing sequence of positive integers
such that $\sum_{i=1}^nd_i=2n-2c$ for some positive integer $c$ and $n_1(d) \leq n_{\geq 2}(d)$,
then there is a forest $F$ with $c$ components and degree sequence $d$ such that
\begin{enumerate}[(i)]
\item there are exactly $n_1(d)$ support vertices in $F$
each of which is adjacent to exactly one vertex of degree $1$,
\item the vertices in $V_{\geq 2}(F)$ that are not support vertices
are all of degree $2$, and induce a path $P$ of order $n-2n_1(d)$, 
\item $\gamma(F)=\left\lceil\frac{n+n_1(d)-2}{3}\right\rceil$,
and $\alpha(F)=\left\lceil\frac{n}{2}\right\rceil$.
\end{enumerate}
\end{lemma}
{\it Proof:} Since $\sum_{i=1}^nd_i=2n-2c$, we obtain
\begin{eqnarray*}
2c 
& = &  2n-\sum_{i=1}^nd_i\\
& = &  -\sum_{i=1}^n(d_i-2)\\
& = &  
-\sum_{i=1}^{n_{\geq 3}(d)}(d_i-2)
-\sum_{i=n_{\geq 3}(d)+1}^{n_{\geq 2}(d)}(d_i-2)
-\sum_{i=n_{\geq 2}(d)+1}^{n}(d_i-2)\\
& = &  -\sum_{i=1}^{n_{\geq 3}(d)}(d_i-2)+n_1(d),
\end{eqnarray*}
which implies
\begin{equation}
n_1(d) = 2c + \sum_{i=1}^{n_{\geq 3}(d)}(d_i-2) \ge 2c + n_{\geq 3}(d).
\label{Eq1}
\end{equation}
By supposition, $n_{\geq 2}(d) \ge n_1(d)$, 
and so, by (\ref{Eq1}), $n_{\geq 2}(d) \ge 2c + n_{\geq 3}(d)$. 
Hence, $n_2(d) = n_{\geq 2}(d) - n_{\geq 3}(d) \ge 2c \ge 2$, implying that $d_{n_{\geq 2}(d)} = 2$. 
By (\ref{Eq1}), $n_1(d) - 2c + 1 > n_{\geq 3}(d)$, and so $d_i \le 2$ for $i \ge n_1(d) - 2c + 1$. 
Therefore, $d_i=2$ for every integer $i$ with $n_1(d)-2c+1\leq i\leq {n_{\geq 2}(d)}$. 
In particular, $d_i=2$ for every integer $i$ with $n_1(d)+1\leq i\leq {n_{\geq 2}(d)}$. 
Thus,
$$d'=(d_1',\ldots,d_{n_1(d)}')=(d_1-1,\ldots,d_{n_1(d)}-1)$$
is a sequence of $n_1(d)$ positive integers with
\begin{eqnarray*}
\sum_{i=1}^{n_1(d)}d_i' 
& = & \left( \sum_{i=1}^{n}d_i - \sum_{i=n_1(d) + 1}^{n_{\geq 2}(d)}d_i - \sum_{i=n_{\geq 2}(d) + 1}^{n}d_i \right) - n_1(d)\\
& = & \sum_{i=1}^{n}d_i-2(n_{\geq 2}(d)-n_1(d))-n_1(d)-n_1(d)\\
&=& \sum_{i=1}^{n}d_i-2n_{\geq 2}(d)\\
&=& 2(n-n_{\geq 2}(d))-2c\\
&=& 2n_1(d)-2c,
\end{eqnarray*}
that is, $d'$ is the degree sequence of a forest $F'$ of order $n_1(d)$ with $c$ components.
Let $F$ arise by
\begin{itemize}
\item attaching one new vertex of degree $1$ to each vertex of $F'$, and
\item subdividing one edge of $F'$ exactly $n-2n_1(d)=n_{\geq 2}(d)-n_1(d)$ times.
\end{itemize}
By construction, $F$ is a forest with degree sequence $d$
that satisfies (i) and (ii).

Since some minimum dominating set of $F$
contains all $n_1(d)$ support vertices of $F$
as well as exactly $\left\lceil\frac{n-2n_1(d)-2}{3}\right\rceil$ interior vertices of $P$,
we obtain
\begin{eqnarray*}
\gamma(F) & = & n_1(d)+\left\lceil\frac{n-2n_1(d)-2}{3}\right\rceil
= \left\lceil\frac{n+n_1(d)-2}{3}\right\rceil.
\end{eqnarray*}
Similarly, some maximum independent set in $F$ contains all vertices of degree $1$
as well as the vertices of the larger partite set of the bipartite graph $P$,
which implies
\begin{eqnarray*}
\alpha(F) & = & n_1(d)+\left\lceil\frac{n-2n_1(d)}{2}\right\rceil
= \left\lceil\frac{n}{2}\right\rceil,
\end{eqnarray*}
and completes the proof. $\Box$.

\medskip

\noindent We proceed to our main results.

Since
\begin{eqnarray*}
\gamma_{\max}^{\cal F}((d_1,\ldots,d_{n-1},0))&=&\gamma_{\max}^{\cal F}((d_1,\ldots,d_{n-1}))+1\mbox{ and}\\
\alpha_{\min}^{\cal F}((d_1,\ldots,d_{n-1},0))&=&\alpha_{\min}^{\cal F}((d_1,\ldots,d_{n-1}))+1,
\end{eqnarray*}
it suffices to consider degree sequences with only positive entries.
Similarly, if $d$ is a degree sequence that contains only $1$-entries,
then $\gamma_{\max}^{\cal F}(d)=\alpha_{\min}^{\cal F}(d)=\frac{n}{2}$.
Therefore, we focus on degree sequences that contain at least one entry that is at least $2$.

\begin{theorem}\label{theorem1}
If $d=(d_1,\ldots,d_n)$ is a non-increasing sequence of positive integers
such that $d_1\geq 2$ and $\sum_{i=1}^nd_i=2n-2c$ for some positive integer $c$, then
$$
\gamma_{\max}^{\cal F}(d)=
\left\{
\begin{array}{rl}
n-n_1(d)+c-1, &
\mbox{if } n_1(d)>n_{\geq 2}(d) \mbox{ and }
c-1<\left\lceil\frac{n_1(d)-n_{\geq 2}(d)}{2}\right\rceil\\[3mm]
\left\lfloor\frac{n}{2}\right\rfloor, &
\mbox{if } n_1(d)>n_{\geq 2}(d)\mbox{ and }
c-1\geq \left\lceil\frac{n_1(d)-n_{\geq 2}(d)}{2}\right\rceil\\[3mm]
\left\lceil\frac{n+n_1(d)-2}{3}\right\rceil, &
\mbox{if } n_1(d)\leq n_{\geq 2}(d).
\end{array}
\right.
$$
\end{theorem}
{\it Proof:} The proof is by induction on $n$.
Since $n+1\leq \sum_{i=1}^nd_i=2n-2c$, we have $n\geq 2c+1$.

If $n=2c+1$, then $d_1=2$ and $d_2=\ldots=d_n=1$.
In this case,
$n$ is odd,
$n_1(d)=n-1>1=n_{\geq 2}(d)$,
and
$c-1 = \frac{n-3}{2}<\frac{n-1}{2} = \left\lceil\frac{n_1(d)-n_{\geq 2}(d)}{2}\right\rceil$.
The only forest with degree sequences $d$ consists of
one component of order $3$ and
$c-1$ components of order $2$,
and has domination number $c=n-n_1(d)+c-1$.

Now, let $n>2c+1$.
As noted above, every forest with degree sequence $d$ has exactly $c$ components.

Let the forest $F$ with degree sequence $d$ be such that
\begin{itemize}
\item $\gamma(F)=\gamma_{\max}^{\cal F}(d)$,
\item subject to the first condition,
the number $k_2$ of components of $F$ of order $2$ is maximum, and
\item subject to the first and second condition,
the number of support vertices of $F$ is maximum.
\end{itemize}

\begin{claim}\label{claim1}
If $c\geq 2$ and $k_2=0$,
then no vertex in $V_{\geq 2}(F)$ has more than one neighbor in $V_1(F)$, and thus $n_1(d)\leq n_{\geq 2}(d)$.
\end{claim}
{\it Proof of Claim \ref{claim1}:}
Suppose that there is a vertex $x$ with two neighbors in $V_1(F)$.
Let $x'$ be a neighbor of $x$ of degree $1$.
Since $c\geq 2$ and $k_2=0$, there is a support vertex $y$ in a component that does not contain $x$.
Let $y'$ be a neighbor of $y$ of degree $1$.
Now, $F'=F-xx'-yy'+xy+x'y'$ is a forest with degree sequence $d$ that has a component of order $2$.
Since $x$ is a support vertex of $F'$,
the forest $F'$ has a minimum dominating set $D'$ with $x,x'\in D'$.
Since $(D'\setminus \{ x'\})\cup \{ y\}$ is a dominating set of $F$,
we obtain
$\gamma_{\max}^{\cal F}(d)\geq \gamma(F')\geq \gamma(F)=\gamma_{\max}^{\cal F}(d)$,
which implies $\gamma(F')=\gamma_{\max}^{\cal F}(d)$,
and yields a contradiction to the choice of $F$.
$\Box$

\begin{claim}\label{claim2}
If $c=1$,
then there are no two vertices $x$ and $y$ in $V_{\geq 2}(F)$
such that $x$ has at least two neighbors in $V_1(F)$ and $y$ has no neighbor in $V_1(F)$.
\end{claim}
{\it Proof of Claim \ref{claim2}:}
Suppose $x$ and $y$ are as in the statement.
Let $x'$ be a neighbor of $x$ of degree $1$,
and let $y'$ be a neighbor of $y$ that does not lie on the path in $F$ between $x$ and $y$.
Now, $F'=F-xx'-yy'+xy'+x'y$ is a forest with degree sequence $d$
that has more support vertices than $F$.
Note that $F'$ is a tree and has a minimum dominating set $D'$ that does not contain a vertex of degree $1$.
Since $D'$ is also a dominating set of $F$,
we obtain $\gamma_{\max}^{\cal F}(d)\geq \gamma(F') \geq \gamma(F) =\gamma_{\max}^{\cal F}(d)$,
which implies $\gamma(F')=\gamma_{\max}^{\cal F}(d)$,
and yields a contradiction to the choice of $F$. $\Box$

\medskip

\noindent We consider two cases.

\medskip

\noindent {\bf Case 1} $n_1(d)>n_{\geq 2}(d)$.

\medskip

\noindent If $c=1$, then
$c-1<\left\lceil\frac{n_1(d)-n_{\geq 2}(d)}{2}\right\rceil$, and, by Claim \ref{claim2}, every vertex in $V_{\geq 2}(F)$ is a support vertex,
which implies $\gamma_{\max}^{\cal F}(d)=\gamma(F)=n_{\geq 2}(d)=n-n_1(d)+c-1$. Hence, we may assume that $c\geq 2$.

Claim \ref{claim1} implies $k_2\geq 1$,
that is, $F$ has a component $K$ of order $2$.
Note that $\gamma(F)=\gamma(F-V(K))+1$,
and that
$F-V(K)$ is a forest with
degree sequence $d'=(d_1,\ldots,d_{n-2})$, $n'=n-2$ vertices, and $c'=c-1$ components.
For $d'$, we obtain
\begin{eqnarray*}
d'_1 &\geq & 2,\\
n_1(d') &=& n_1(d)-2,\\
n_{\geq 2}(d')&=&n_{\geq 2}(d),\mbox{ and}\\
\sum_{i=1}^{n'}d_i &=& 2n'-2c'.
\end{eqnarray*}
By the choice of $F$,
we have $\gamma(F-V(K))=\gamma_{\max}^{\cal F}(d')$,
which implies
$\gamma_{\max}^{\cal F}(d)=\gamma_{\max}^{\cal F}(d')+1$.

First, we assume that $n_1(d')>n_{\geq 2}(d')$.

If $c-1\geq \left\lceil\frac{n_1(d)-n_{\geq 2}(d)}{2}\right\rceil$,
then $c'-1\geq \left\lceil\frac{n_1(d')-n_{\geq 2}(d')}{2}\right\rceil$,
and, by induction,
$$\gamma_{\max}^{\cal F}(d)
=\left\lfloor\frac{n'}{2}\right\rfloor +1
=\left\lfloor\frac{n}{2}\right\rfloor.$$

If $c-1<\left\lceil\frac{n_1(d)-n_{\geq 2}(d)}{2}\right\rceil$,
then $c'-1<\left\lceil\frac{n_1(d')-n_{\geq 2}(d')}{2}\right\rceil$, and, by induction,
$$\gamma_{\max}^{\cal F}(d)
=n'-n_1(d')+c'-1+1
=n-n_1(d)+c-1.$$

Next, we assume that $n_1(d')\leq n_{\geq 2}(d')$.

In this case, $n_1(d)>n_{\geq 2}(d)$
implies that $n_1(d)\in \{ n_{\geq 2}(d)+1,n_{\geq 2}(d)+2\}$.
This implies $c-1\geq 1=\left\lceil\frac{n_1(d)-n_{\geq 2}(d)}{2}\right\rceil$, and $n=2n_1(d)-r$ for some $r\in \{ 1,2\}$.

By induction, we obtain
\begin{eqnarray*}
\gamma_{\max}^{\cal F}(d) &=&\gamma_{\max}^{\cal F}(d')+1\\
&=& \left\lceil\frac{n'+n_1(d')-2}{3}\right\rceil+1\\
&=& \left\lceil\frac{n+n_1(d)-6}{3}\right\rceil+1\\
&=& \left\lceil\frac{3n_1-r-3}{3}\right\rceil\\
&=& n_1-1\\
&=& \left\lfloor\frac{n}{2}\right\rfloor.
\end{eqnarray*}
Altogether, in each case,
$\gamma_{\max}^{\cal F}(d)$ has the value stated in the theorem.

\medskip

\noindent {\bf Case 2} $n_1(d)\leq n_{\geq 2}(d)$.

\medskip

\noindent By Lemma \ref{lemma2},
$\gamma_{\max}^{\cal F}(d)\geq \left\lceil\frac{n+n_1(d)-2}{3}\right\rceil$, and it remains to show
$\gamma_{\max}^{\cal F}(d)\leq \left\lceil\frac{n+n_1(d)-2}{3}\right\rceil$.

First, we assume that $k_2=0$.
By Claim \ref{claim1} and Claim \ref{claim2},
we obtain that no vertex in $V_{\geq 2}(F)$
has more than one neighbor in $V_1(F)$.
Therefore,
if $U$ is the set of vertices in $V_{\geq 2}(F)$
that are not support vertices,
then $|U|=n_{\geq 2}(d)-n_1(d)$.
Let the subgraph $F'$ of $F$ induced by $U$
have components of orders $p_1,\ldots,p_k$, respectively.
Note that every vertex of degree at most $1$ in $F'$
is adjacent to a support vertex of $F$.
Therefore, the set of the $n_1(d)$ support vertices of $F$
together with sets as in Lemma \ref{lemma1} for each component of $F'$
form a dominating set of $F$.
By Lemma \ref{lemma1}, we obtain
\begin{eqnarray*}
\gamma_{\max}^{\cal F}(d) &=& \gamma(F)\\
& \leq & n_1(d)+
\left\lceil\frac{p_1-2}{3}\right\rceil
+\cdots+
\left\lceil\frac{p_k-2}{3}\right\rceil\\
& \leq & n_1(d)+\left\lceil\frac{(p_1+\cdots+p_k)-2}{3}\right\rceil\\
& = & n_1(d)+\left\lceil\frac{n_{\geq 2}(d)-n_1(d)-2}{3}\right\rceil\\
& = & \left\lceil\frac{n+n_1(d)-2}{3}\right\rceil.
\end{eqnarray*}
Next, we assume that $k_2\geq 1$.
As in Case 1, this implies
$\gamma_{\max}^{\cal F}(d)=\gamma_{\max}^{\cal F}((d_1,\ldots,d_{n-2}))+1$.
Since $n_1(d)-2<n_{\geq 2}(d)$,
we obtain, by induction,
\begin{eqnarray*}
\gamma_{\max}^{\cal F}(d) & = & \gamma_{\max}^{\cal F}((d_1,\ldots,d_{n-2}))+1\\
& = & \left\lceil\frac{(n-2)+(n_1(d)-2)-2}{3}\right\rceil+1\\
& = & \left\lceil\frac{n+n_1(d)-3}{3}\right\rceil\\
& \leq & \left\lceil\frac{n+n_1(d)-2}{3}\right\rceil,
\end{eqnarray*}
which completes the proof. $\Box$

\begin{theorem}\label{theorem2}
If $d=(d_1,\ldots,d_n)$ is a non-increasing sequence of positive integers
such that $d_1\geq 2$ and $\sum_{i=1}^nd_i=2n-2c$ for some positive integer $c$, then
$$
\alpha_{\min}^{\cal F}(d)=
\left\{
\begin{array}{rl}
n_1(d)-c+1, & \mbox{if }
n_1(d)>n_{\geq 2}(d)\mbox{ and }
c-1<\left\lceil\frac{n_1(d)-n_{\geq 2}(d)}{2}\right\rceil\\[3mm]
\left\lceil\frac{n}{2}\right\rceil, & \mbox{ otherwise.}
\end{array}
\right.
$$
\end{theorem}
{\it Proof:} The proof is by induction on $n$ and quite similar to the proof of Theorem \ref{theorem1}.

Again $n\geq 2c+1$, and if $n=2c+1$,
then
$n_1(d)>n_{\geq 2}(d)$,
$c-1<\left\lceil\frac{n_1(d)-n_{\geq 2}(d)}{2}\right\rceil$,
and the unique forest with degree sequences $d$ has independence number $c+1=n_1(d)-c+1$.

Now, let $n>2c+1$.

Let the forest $F$ with degree sequence $d$ be such that
\begin{itemize}
\item $\alpha(F)=\alpha_{\min}^{\cal F}(d)$,
\item subject to the first condition,
the number $k_2$ of components of $F$ of order $2$ is maximum, and
\item subject to the first and second condition,
the number of support vertices of $F$ is maximum.
\end{itemize}

\begin{claim}\label{claim3}
If $c\geq 2$ and $k_2=0$,
then no vertex in $V_{\geq 2}(F)$ has more than one neighbor in $V_1(F)$, and thus $n_1(d)\leq n_{\geq 2}(d)$.
\end{claim}
{\it Proof of Claim \ref{claim3}:}
Let $x$, $x'$, $y$, $y'$, and $F'$ be exactly as in the proof of Claim \ref{claim1}.
The forest $F'$ has a maximum independent set $I'$ that contains $V_1(F)\cap N_F(x)$. 
We note that $I'$ contains $x'$ but contains neither $x$ nor $y'$. 
The set $I'$ possibly contains $y$.
Now, $(I'\setminus \{ y\})\cup \{ y'\}$ is an independent set of $F$,
which implies
$\alpha_{\min}^{\cal F}(d)\leq \alpha(F')\leq \alpha(F)= \alpha_{\min}^{\cal F}(d)$,
and so, $\alpha(F')=\alpha_{\min}^{\cal F}(d)$.
Since $F'$ has a component of order $2$, we obtain a contradiction to the choice of $F$.
$\Box$

\begin{claim}\label{claim4}
If $c=1$,
then there are no two vertices $x$ and $y$ in $V_{\geq 2}(F)$
such that $x$ has at least two neighbors in $V_1(F)$ and $y$ has no neighbor in $V_1(F)$.
\end{claim}
{\it Proof of Claim \ref{claim4}:}
Let $x$, $x'$, $y$, $y'$, and $F'$ be exactly as in the proof of Claim \ref{claim2}.
Some maximum independent set $I'$ in $F'$ contains $x'$ and a neighbor of degree $1$ of $x$.
Since $I'$ is independent in $F$, we obtain $\alpha(F')=\alpha_{\min}^{\cal F}(d)$.
Since $F'$ has more support vertices than $F$, we obtain a contradiction to the choice of $F$. $\Box$

\medskip

\noindent We consider two cases.

\medskip

\noindent {\bf Case 1} $n_1(d)>n_{\geq 2}(d)$.

\medskip

\noindent If $c=1$, then
$c-1<\left\lceil\frac{n_1(d)-n_{\geq 2}(d)}{2}\right\rceil$, and, by Claim \ref{claim4}, every vertex in $V_{\geq 2}(F)$ is a support vertex,
which implies $\alpha_{\min}^{\cal F}(d)=\alpha(F)= n_1(d)-c+1$. Hence, we may assume that $c\geq 2$.

Claim \ref{claim3} implies $k_2\geq 1$,
and hence, $\alpha_{\min}^{\cal F}(d)=\alpha_{\min}^{\cal F}(d')+1$
where $d'=(d_1,\ldots,d_{n-2})$ has the same properties as stated in the proof of Theorem \ref{theorem1}.

First, we assume that $n_1(d')>n_{\geq 2}(d')$.

If $c-1\geq \left\lceil\frac{n_1(d)-n_{\geq 2}(d)}{2}\right\rceil$,
then $c'-1\geq \left\lceil\frac{n_1(d')-n_{\geq 2}(d')}{2}\right\rceil$,
and, by induction,
$\alpha_{\min}^{\cal F}(d)=\left\lceil\frac{n'}{2}\right\rceil +1=\left\lceil\frac{n}{2}\right\rceil$.

If $c-1<\left\lceil\frac{n_1(d)-n_{\geq 2}(d)}{2}\right\rceil$,
then $c'-1<\left\lceil\frac{n_1(d')-n_{\geq 2}(d')}{2}\right\rceil$, and, by induction,
$\alpha_{\min}^{\cal F}(d)=n_1(d')-c'+1+1=n_1-c+1$.

Next, we assume that $n_1(d')\leq n_{\geq 2}(d')$.
In this case, $c-1\geq 1=\left\lceil\frac{n_1(d)-n_{\geq 2}(d)}{2}\right\rceil$.
By induction, $\alpha_{\min}^{\cal F}(d)=\left\lceil\frac{n'}{2}\right\rceil +1=\left\lceil\frac{n}{2}\right\rceil$.
Altogether, in each case,
$\alpha_{\min}^{\cal F}(d)$ has the value stated in the theorem.

\medskip

\noindent {\bf Case 2} $n_1(d)\leq n_{\geq 2}(d)$.

\medskip

\noindent By Lemma \ref{lemma2},
$\alpha_{\min}^{\cal F}(d)\leq \left\lceil\frac{n}{2}\right\rceil$, and it remains to show
$\alpha_{\min}^{\cal F}(d)\geq \left\lceil\frac{n}{2}\right\rceil$.

First, we assume that $k_2=0$.
By Claim \ref{claim3} and Claim \ref{claim4},
we obtain that no vertex in $V_{\geq 2}(F)$
has more than one neighbor in $V_1(F)$.
Let $U$ and $F'$ be exactly as in the proof of Theorem \ref{theorem1}.
The set of the $n_1(d)$ vertices of degree $1$ of $F$
together with the larger partite set of the bipartite graph $F'$
form an independent set in $F$, and we obtain
\begin{eqnarray*}
\alpha_{\min}^{\cal F}(d) &=& \alpha(F)\\
& \geq & n_1(d)+\left\lceil\frac{|U|}{2}\right\rceil\\
& = & n_1(d)+\left\lceil\frac{n_{\geq 2}(d)-n_1(d)}{2}\right\rceil\\
& = & \left\lceil\frac{n}{2}\right\rceil.
\end{eqnarray*}
Next, we assume that $k_2\geq 1$.
As in Case 1, this implies
$\alpha_{\min}^{\cal F}(d)=\alpha_{\min}^{\cal F}((d_1,\ldots,d_{n-2}))+1$.
Since $n_1(d)-2<n_{\geq 2}(d)$,
we obtain, by induction,
$\alpha_{\min}^{\cal F}(d)=\left\lceil\frac{n-2}{2}\right\rceil+1=\left\lceil\frac{n}{2}\right\rceil$,
which completes the proof. $\Box$

\medskip
\noindent It is a curious fact, that
$\gamma_{\max}^{\cal F}(d)+\alpha_{\min}^{\cal F}(d)=n$ for $n_1(d)>n_{\geq 2}(d)$.

\end{document}